\documentclass[11pt, reqno]{amsart}

\usepackage{amssymb,latexsym,amsmath,amsfonts}
\usepackage{mathrsfs}
\usepackage{graphicx}

\numberwithin{equation}{section}

\theoremstyle{definition}
\newtheorem{definition}{Definition}[section]

\theoremstyle{remark}
\newtheorem{remark}[definition]{Remark}

\theoremstyle{plain}
\newtheorem{theorem}[definition]{Theorem}
\newtheorem{result}[definition]{Result}
\newtheorem{lemma}[definition]{Lemma}
\newtheorem{proposition}[definition]{Proposition}

\newtheorem{key}[definition]{Key Lemma}

\newcommand{\eps}{\varepsilon}
\newcommand{\zt}{\zeta}
\newcommand{\zbar}{\overline{z}}

\newcommand{\bas}{\boldsymbol{\epsilon}}

\newcommand{\bdy}{\partial}
\newcommand{\OM}{\Omega}

\newcommand{\ome}{\omega}

\newcommand{\smoo}{\mathcal{C}}
\newcommand{\hol}{\mathcal{O}}

\newcommand\ber[1]{\boldsymbol{{\sf B}}_{{#1}}}
\newcommand\con[1]{{\sf const}_{{#1}}}
\newcommand{\dee}{\boldsymbol{{\sf D}}}
\newcommand{\qu}{\boldsymbol{{\sf Q}}}
\newcommand{\eff}{\boldsymbol{{\sf f}}}

\newcommand{\re}{{\sf Re}}

\newcommand{\bcdot}{\boldsymbol{\cdot}}

\newcommand{\lrarw}{\longrightarrow}

\newcommand{\cjac}{{\sf Jac}_{\mathbb{C}}}

\newcommand{\bund}{\mathfrak{B}}
\newcommand{\str}{\mathcal{M}}

\newcommand{\Cn}{\mathbb{C}^n}
\newcommand{\C}{\mathbb{C}} 
\newcommand{\R}{\mathbb{R}}
\newcommand{\Z}{\mathbb{Z}}
\newcommand{\N}{\mathbb{N}}
\newcommand{\pro}{\mathbb{C}\mathbb{P}}

\begin{document}

\title[Proper holomorphic maps]{Proper holomorphic maps between bounded symmetric domains revisited}

\author{Gautam Bharali}
\address{Department of Mathematics, Indian Institute of Science, Bangalore 560012, India}
\email{bharali@math.iisc.ernet.in}

\author{Jaikrishnan Janardhanan}
\address{Department of Mathematics, Indian Institute of Science, Bangalore 560012, India}
\email{jaikrishnan@math.iisc.ernet.in}

\thanks{GB is supported by a UGC Centre of Advanced Study grant. JJ is supported by a UGC Centre for
Advanced Study grant and by a scholarship from the IISc}


\begin{abstract}
We prove that a proper holomorphic map between two non-planar bounded symmetric domains
of the same dimension, one of them being irreducible, is a biholomorphism. Our methods allow us
to give a single, all-encompassing argument that unifies the various special cases in which
this result is known. We discuss an application of these methods to domains having
noncompact automorphism groups that are not assumed to act transitively.
 
\end{abstract}
\maketitle

\section{Introduction and statement of results}\label{S:intro}

The primary objective of this paper is to prove the following result:

\begin{theorem}\label{T:mThm}
 Let $D_1$ and $D_2$ be two bounded symmetric domains of complex dimension $n\geq 2$.
 Assume that either $D_1$ or $D_2$ is irreducible. Then, any proper holomorphic mapping of $D_1$ into
 $D_2$ is a biholomorphism.
\end{theorem}

The above theorem is known in several special cases. For $D_1 = D_2 = \mathbb{B}^n$,
the (Euclidean) ball in $\Cn, \ n\geq 2$, the result was established by Alexander 
\cite{alexander1977:phm}. This is a pioneering work that has motivated several generalizations to
proper holomorphic maps between certain types of smoothly bounded pseudoconvex domains. 
Henkin and Novikov \cite{henkinNovikov1984:pmcd} described a method for proving the above result
when $D_1 = D_2$\,($=D$, say) and $D$ is a bounded symmetric domain that is {\em not}
of tube type. About a decade later, Tsai \cite{tsai1993:rigidity} established the result for
$D_1$ and $D_2$ as above, provided $D_1$ is irreducible
and ${\rm rank}(D_1)\geq {\rm rank}(D_2)\geq 2$.
\smallskip

Tsai's result is a broad metric-rigidity theorem (under the Bergman metric) for proper holomorphic maps
of $D_1$ into $D_2$, where $D_1$ and $D_2$ are as above but not necessarily
equidimensional. In such a result, the condition ${\rm rank}(D_1)\geq {\rm rank}(D_2)\geq 2$
is indispensible. Adapting Tsai's ideas to the equidimensional
case, Tu \cite{tu2002:rigidity} established Theorem~\ref{T:mThm} in the higher-rank case, assuming
$D_1$ {\em is irreducible}. In using Tsai's ideas, the assumption that $D_1$ is irreducible is
essential --- see \cite[Proposition\,3.3]{tu2002:rigidity} --- and it is not clear that a small mutation
of those ideas allows one to weaken this assumption. In our work, we are able to assume
either $D_1$ or $D_2$ to be irreducible precisely by not relying too heavily on the fine structure of
these domains.  Indeed, we wish to emphasize that the focus of this work is
{\em not} on mopping up the residual cases in Theorem~\ref{T:mThm}. The methods in
\cite{henkinNovikov1984:pmcd} (and \cite{tumanovKhenkin82:localChar}, on
which \cite{henkinNovikov1984:pmcd} relies), \cite{tsai1993:rigidity} and  \cite{tu2002:rigidity}
are tied, in a rather maximalistic way, to the fine structure of a bounded symmetric domain. In
contrast, we present some ideas that make very mild use of the underlying orbit structure
of the bounded symmetric domains. They could therefore be applied to manifolds
whose automorphism groups are not assumed to act transitively but are merely ``large enough''.
Theorem~\ref{T:short} is an illustration of this notion. These ideas also provide a {\em unified}
argument, irrespective of rank or reducibility, for Theorem~\ref{T:mThm}.    
\smallskip

We need to be more precise about the preceding remarks. This requires some elaboration on
the objects of interest. A bounded symmetric domain in $\Cn$ is the holomorphic imbedding in
$\Cn$ of some Hermitian symmetric space of noncompact type. It is {\em irreducible} if
it is not a product of bounded symmetric domains of lower dimension. Cartan studied
Hermitian symmetric spaces of noncompact type and classified the irreducible ones, showing
that they are one of six types of homogeneous spaces. An outcome of Harish-Chandra's
work in \cite{harishchandra1956:sslg-VI} is that these homogeneous spaces (and products
thereof) can be imbedded in $\Cn$ as
bounded convex balanced domains (we say that a domain $D\subset \Cn$
is {\em balanced} if, for any $z\in D$, $\zt z\in D$ for each $\zt$ in the closed unit disc
centered at $0\in \C$). This imbedding is unique up to a linear isomorphism of $\Cn$. Such
a realization of a bounded symmetric domain is called a {\em Harish-Chandra realization}.
\smallskip

The three main features that we wish to emphasize about this work are:
\begin{itemize}
 \item[$a)$] The arguments in \cite{alexander1977:phm} involve many estimates showing
 how a proper mapping maps conical regions with vertex on $\bdy\mathbb{B}^n$
 into the ``admissible'' approach regions of Kor{\'a}nyi
 \cite{koranyi1969:hhh}. Boundary approach, in a somewhat different sense, plus
 Chern--Moser theory \cite{chernMoser1974:realHypersurfaces} make an appearance
 in \cite{tumanovKhenkin82:localChar}. In contrast, apart from, and owing to, a result
 of Bell \cite{bell1982:proper-circular} on boundary behaviour, our proof
 involves rather ``soft'' methods.
 
  \item[$b)$] The techniques underlying \cite{tsai1993:rigidity} and \cite{tu2002:rigidity}
  rely almost entirely on the fine structure of a bounded symmetric domain. Specifically,
  they involve studying the effect of a proper holomorphic map on the characteristic
  symmetric subspaces of a bounded symmetric domain of rank\,$\geq 2$. In contrast,
  our techniques rely on only a coarse distinction between the different strata that
  comprise the boundary of an irreducible bounded symmetric domain (e.g.,
  see Remark~\ref{rem:maxStratum} below).
  
  \item[$c)$] An advantage of arguments that rely on only a coarse resolution of 
  a bounded symmetric domain is that some of them are potentially applicable to the study
  of domains that have noncompact automorphism groups, but are not assumed to be
  symmetric. A demonstration this viewpoint is the proof of Theorem~\ref{T:short} below.
\end{itemize}

Let $D$ be a bounded symmetric domain. The main technical tool that facilitates our
study of the structure of $\bdy D$, and describes certain elements of ${\rm Aut}(D)$
with the optimal degree of explicitness, is the notion of Jordan triple systems. The application
of Jordan triple systems to geometry appears to have been pioneered by Koecher
\cite{koecher1999:minnesota}. Our reference on this subject are the lecture notes of
Loos \cite{loos1977:bsdJp}, which are devoted specifically to the bounded symmetric domains.
Jordan triple systems and versions of the Schwarz lemma are our primary tools. We
present next an outline of how we use these tools.
\smallskip

An important lemma, which is inspired by Alexander's work, is the following

\begin{key}
 Let $D$ be a realization of an irreducible bounded symmetric domain of dimension $n\geq 2$
 as a bounded convex balanced domain  in $\Cn$.
 For $z\in D\setminus\{0\}$, let $\Delta_z:=\{\zt z: \zt\in \C \ \text{and} \ \zt z\in D\}$.
 Let $W_1$ and $W_2$ be two regions in $D$ such that $0\in W_1\cap W_2$ and let
 $F : D\to D$ be a holomorphic map. Assume that:
 \begin{itemize}
  \item[$i)$]  $F$ maps $W_1$ biholomorphically onto $W_2$ with $F(0)=0$.
  \item[$ii)$] There exists a non-empty open set $U\subset W_1\setminus \{0\}$ such that,
  for each $z\in U$, $\Delta_z\subset W_1$ and $\Delta_{F(z)}\subset W_2$.
 \end{itemize}
 Then, $F$ is an automorphism of $D$.
\end{key}

\noindent This is a consequence of Vigu{\'e}'s Schwarz lemma
\cite{vigue1991:domainesSymetriques} (see Result~\ref{R:SL-vig} below), and the irreducibility
of $D$ is essential to the lemma.
\smallskip 

Our proof of Theorem~\ref{T:mThm} may be summarized as follows (we will
assume here that $D_1$ and $D_2$ are Harish-Chandra realizations of the domains in question): 
\begin{itemize}
 \item By Bell's theorem \cite[Theorem~2]{bell1982:proper-circular}, $F$ extends to
 a neighbourhood of $\overline{D}_1$ and we can find a point
 $p$ in the Bergman--Shilov boundary of $D_1$, and a small ball $B$ around it,
 such that $F|_B$ is a biholomorphism.
 \item We may assume that $F(0)=0$.
 Let $\{a_k\}$ be a sequence in $D_1\cap B$ converging to $p$ and let
 $b_k:=F(a_k)$. Let $\phi^j_k\in {\rm Aut}(D_j)$ be an automorphism that maps
 $0$ to $a_k$ if $j=1$, and to $b_k$ if $j=2$. It turns out that both $p$ and $F(p)$
 are peak points, whence $\phi^j_k\lrarw p^{(j)}$ uniformly on compact subsets,
 where $p^{(1)}:=p$ and $p^{(2)}:=F(p)$.
 \item Using the Schwarz lemma for convex balanced domains (Result~\ref{R:SL-cbd}
 below) we show that a subsequence of $\{(\phi^2_j)^{-1}\circ F\circ \phi^1_j\}$
 converges to a linear map and that, owing to the tautness of $D_1$ and $D_2$, this
 map is a biholomorphism of $D_1$ onto $D_2$.
 \item We may now take  $D_1=D_2=D$. We shall use our Key Lemma, with
 $W_1=(\phi^1_k)^{-1}(D\cap B)$ and $W_2=(\phi^2_k)^{-1}(D\cap F(B))$ for
 $k$ sufficiently large. 
 \item Since the analytic discs $\Delta_z$ and $\Delta_{F(z)}$ are not relatively
 compact in $D$, the mode of convergence of $\{\phi^j_k\}$
 isn't {\em a priori} good enough to infer that appropriate families of these
 discs will be swallowed up by $W_j$, $j=1,2$. By Bell's theorem, each 
 $\phi^j_k$ extends to some neighbourhood of $D$. We show that
 $\{\phi^j_k\}$, passing to a subsequence and relabelling if necessary, converges
 {\em uniformly} on certain special special circular subsets of $D$ that are
 adherent to $\bdy D$. This is enough to overcome the difficulty just described.
\end{itemize} 

Let us define a term that we used in the sketch above, which we shall also need
in stating our next theorem.

\begin{definition}
 Let $D\varsubsetneq \Cn$ be a domain and let $p\in \bdy D$. We say that
 $p$ is a {\em peak point} if there exists a function 
 $h\in \hol(D)\cap\smoo(\overline{D};\C)$ such that $h(p)=1$ and
 $|h(z)| < 1 \ \forall z\in \overline{D}\setminus\{p\}$. The function 
 $h$ is called a {\em peak function for $p$}.
\end{definition}

When a domain $D$ is bounded, the noncompactness of ${\rm Aut}(D)$ (in the compact-open
topology) is equivalent to $D$ having a boundary orbit-accumulation point;
see \cite{narasimhan1971:SCV}.

\begin{definition}
 Let $D\varsubsetneq \Cn$ be a domain and let $p\in \bdy D$. We say that
 $p$ is a {\em boundary orbit-accumulation point} if there exist a point
 $a\in D$ and a sequence of automorphisms $\{\phi_k\}$ of $D$ such that
 $\lim_{k\to \infty}\phi_k(a)=p$.
\end{definition}

With the last two definitions, we are in a position to state our second theorem.
Note that $D_1$ is {\em not} assumed to be a bounded symmetric domain. Yet, some
of the techniques sketched above (versions of which have been used to remarkable
effect in the literature in this field) are general enough to be applicable to the following
situation.

\begin{theorem}\label{T:short}
 Let $D_1$ be a bounded convex balanced domain in $\Cn$ whose
 automorphism group is noncompact and let $p$ be a boundary orbit-accumulation
 point. Let $D_2$ be a realization of
 a bounded symmetric domain as a bounded convex balanced domain in $\Cn$.
 Assume that there is a neighbourhood $U$ of $p$ and a biholomorphic
 map $F: U\to \Cn$ such that $F(U\cap D_1)\subset D_2$ and
 $F(U\cap \partial{D_1})\subset \partial{D_2}$. Assume that either $p$ or
 $F(p)$ is a peak point. Then, there exists a {\em linear} map that maps
 $D_1$ biholomorphically onto $D_2$.
\end{theorem}

 \begin{remark}\label{rem:rigid}
 Theorem~\ref{T:short} (together with Bell's theorem \cite{bell1982:proper-circular})
 gives a very short proof of the rigidity theorem of 
 Mok and Tsai \cite{mokTsai1992:convexRigidity} under the additional assumption that the
 convex domain $D$ in their result is also circular. There is an extensive literature on
 rigidity theorems relating to bounded symmetric domains, but we shall not dwell any further on it.
\end{remark}

\begin{remark}
 We remark that a version of the above result can be proved without assuming that
 $D_1$ is either balanced or convex. $D_1$ merely needs to be complete Kobayashi
 hyperbolic. However, in this case, the biholomorphism of $D_1$ onto $D_2$
 will not, in general, be linear. We prefer the above version: the conclusion that
 there exists a {\em linear} equivalence places Theorem~\ref{T:short} among the
 rigidity theorems alluded to in Remark~\ref{rem:rigid}. 
\end{remark}
\smallskip

The layout of this paper is as follows. Since Jordan triple systems play a vital role
in describing not just the structure of the boundary of a bounded symmetric domain,
but also some of its key automorphisms, we begin with a primer
on Jordan triple systems. Readers who are familiar with Jordan triple systems
can skip to Section~\ref{S:bdry}, where we discuss the boundary geometry of
bounded symmetric domains. Section~\ref{S:essn} is devoted to stating
and proving certain propositions that are essential to our proofs. Finally, in
Sections \ref{S:mThm} and \ref{S:short}, we present the proofs of the results
stated above
\medskip

\section{A primer on Jordan triple systems}\label{S:primer}

There is  a natural connection between bounded symmetric domains and
certain Hermitian Jordan triple systems. This section collects several
definitions and results that are required to give a coherent description of the
boundary of a bounded symmetric domain (which we shall discuss in
the next section).
\smallskip

Unless otherwise stated, the results in this section can be found in the
UC-Irvine lectures by Loos \cite{loos1977:bsdJp} describing how Jordan triple systems can
be used to study the geometry of bounded symmetric domains.

\begin{definition}
  A \emph{Hermitian Jordan triple system} is a complex vector space $V$ endowed
  with a triple product $(x,y,z) \longmapsto \{x,y,z\}$ that is symmetric and
  bilinear in $x$ and $z$ and conjugate-linear in $y$, and satisfies the Jordan
  identity
  \begin{align*}
	\{x,y,\{u,v,w\}\} &- \{u,v,\{x,y,w\}\}\\
	        &= \{\{x,y,u\},v,w\} - \{u,\{y,x,v\},w\} \ \ \forall x,y,u,v,w \in V.
  \end{align*}
  Such a system is said to be {\em positive} if for each $x \in V \setminus \{0\}$ for
  which $\{x,x,x\} = \lambda x$ (where $\lambda$ is a scalar), we have $\lambda >
  0$.
\end{definition}

As mentioned in Section~\ref{S:intro}, a bounded symmetric domain of complex
dimension $n$ has a realization $D$ as a  bounded convex balanced domain in
$\Cn$. Let $(z,\dots,z_n)$ be the global holomorphic coordinates coming from the
product structure on $\Cn$ and let $(\bas_1,\dots,\bas_n)$ denote the standard
ordered basis of $\Cn$. Let $K_D$ denote the Bergman kernel of (the above
realization of) $D$ and $h_D$ the Bergman metric. The function
$\{\bcdot,\bcdot,\bcdot\} : \Cn \times \Cn \times \Cn \to \Cn$ obtained by
the requirement 
\begin{equation}\label{E:corrJTS}
  h_D(\{\bas_i,\bas_j,\bas_k\}, \bas_l) = \left.\frac{\partial^4 \log
  K_D(z,z)}{\partial z_i\partial\zbar_j\partial z_k \partial \zbar_l}\right|_{z=0},
\end{equation}
and by extending $\C$-linearly in the first and third variables and
$\C$-antilinearly in the second, has the property that
$(\Cn,\{\bcdot,\bcdot,\bcdot\})$ is a positive Hermitian Jordan triple system
(abbreviated hereafter as PHJTS). This relationship is a one-to-one correspondence
between finite-dimensional PHJTSs and bounded symmetric domains --- which we shall
make more precise in Section~\ref{S:bdry}.
\smallskip

Let $(V,\{\bcdot,\bcdot,\bcdot\})$ be a HJTS. It will be convenient to work with
the operators
\begin{equation}\label{E:ops}
  \dee(x,y)z = \qu(x,z)y := \{x,y,z\}.
\end{equation}
We define the operator $Q:V \to \text{End}(V)$ by $Q(x)y := \qu(x,x)y/2$. For any $x \in
V$, we can define the so-called \emph{odd powers} of $x$ recursively by:
\[
  x^{(1)} := x \quad\text{and} \quad x^{(2p+1)} := Q(x)x^{(2p-1)} \text{ if } p \geq 1.
\]
A vector $e \in V$ is called a tripotent if $e^{(3)} = e$.
\smallskip

Tripotents are important to this discussion because:
\begin{itemize}
  \item A finite-dimensional PHJTS has plenty of non-zero tripotents.
  \item Given a finite-dimensional PHJTS $(V,\{\bcdot,\bcdot,\bcdot\})$, any
	vector $V$ has a certain canonical decomposition as a linear combination of
	tripotents.
  \item In a finite-dimensional PHJTS, the set of tripotents forms a
	real-analytic submanifold.
\end{itemize}
We refer the interested reader to \cite[Chapter 3]{loos1977:bsdJp} for details of the first
fact. As for the second fact, we need a couple of new notions. First:
given a HJTS $(V,\{\bcdot,\bcdot,\bcdot\})$, we say that two tripotents $e_1,e_2
\in V$ are orthogonal if $\dee(e_1,e_2) = 0$. Second: given $x \in V$, we define
the real vector space $\ll\!x\!\gg$ by 
\[
  \ll\!x\!\gg\,:=\,\text{span}_{\R}\{ x^{(2p+1)} : p = 0,1,2,\dots\}.
\]
These two notions allows us to state the following:
\begin{result}[Spectral decomposition theorem]\label{R:spec} Let
$(V,\{\bcdot,\bcdot,\bcdot\})$ be a finite-dimensional PHJTS. Then, each $x
\in V\setminus\{0\}$ can be written uniquely as 
\begin{equation}\label{E:spec}
  x = \lambda_1 e_1 + \dots + \lambda_s e_s
\end{equation}
where $\lambda_1 > \lambda_2 > \dots > \lambda_s > 0$ and $\{e_1,\dots,e_s\}$ is
a $\R$-basis of $\ll\!x\!\gg$ comprising pairwise orthogonal tripotents.  
\end{result}

The decomposition of $x \in V$ as given by Result~\ref{R:spec} is called the
\emph{spectral decomposition} of $x$. The assignment $x \longmapsto \lambda_1(x)$,
where $\lambda_1(x)$ is as given by \eqref{E:spec}, is a
well-defined function and can be shown to be a norm on $V$. This norm is called the
\emph{spectral norm} on $V$.
\smallskip

Next, we present another decomposition, which give us the second
ingredient needed to describe the boundary geometry of a bounded symmetric domain.

\begin{result}[Pierce decomposition]\label{R:pierce} Let $(V,\{\bcdot,\bcdot,\bcdot\})
  $ be a HJTS and let $e \in V$ be a tripotent. Then, the spectrum of $\dee (e,e)$ 
  is a subset of $\{0,1,2\}$. Let 
  \[
	V_j = V_j(e) := \{ x \in V : \dee (e,e)x = jx\}, \; \; j \in \mathbb{Z}.
  \]
  Then:
  \begin{enumerate}
	\item[$a)$] $V = V_0 \oplus V_1 \oplus V_2$.
	\item[$b)$] If $e \neq 0$, then $e \in V_2$.
	\item[$c)$] We have the relation
	  $\{V_{\alpha},V_{\beta},V_{\gamma}\} \subset V_{\alpha-\beta + \gamma}$.
	\item[$d)$] $V_0,V_1$ and $V_2$ are Hermitian Jordan subsystems of
	  $\{\bcdot,\bcdot,\bcdot\}$.
  \end{enumerate}
\end{result}
The direct-sum decomposition $(a)$ given by the above result is called the
\emph{Pierce decomposition of $V$ with respect to the tripotent $e$}. The ideas
that go into proving the Pierce decomposition theorem allow us to construct a
special partial order on the set of tripotents of $V$. In order to avoid
statements that are vacuously true, unless stated otherwise, we take
$(V,\{\bcdot,\bcdot,\bcdot\})$ to be a PHJTS. Let $e,e' \in V$ be tripotents. We say 
that \emph{$e$ is dominated by $e'$} ($e \preceq e'$) if there is a tripotent $e_1$
orthogonal to $e$ such that $e' = e + e_1$. We say that \emph{$e$ is strongly 
dominated by $e'$} ($e \prec e'$) if $e \preceq e'$ and $e \neq e'$. The result of 
interest, in this regard, is the following:

\begin{result}\label{R:po}
  Let $(V,\{\bcdot,\bcdot,\bcdot\})$ be a HJTS. Let $e_1,e_2 \in V$ be
  orthogonal tripotents and let $e = e_1 + e_2$. If $e' \in V$ is a tripotent
  orthogonal to $e$, then $e'$ is orthogonal to $e_1$ and $e_2$.
  \smallskip
  
  Now suppose $\{\bcdot,\bcdot,\bcdot\}$ is positive. Then, the relation
  $\preceq$ is a partial order on the set of tripotents.
\end{result}

\begin{definition}
  A tripotent is said to be \emph{minimal} (or \emph{primitive}) if it is minimal for
  $\preceq$ among non-zero tripotents. It is said to be \emph{maximal} if it is 
  maximal for $\preceq$.
\end{definition}

\begin{result}\label{R:maxmin}
  Consider the tripotents of $V$ partially ordered by $\preceq$.
  \begin{enumerate}
	\item A tripotent $e$ is maximal if and only if the Pierce space $V_0(e) =
	  0$.
	\item If, for a tripotent $e$, the Pierce space $V_2(e) = \C e$, then $e$ is
	  primitive.
  \end{enumerate}
\end{result}

Let us now also assume that $(V,\{\bcdot,\bcdot,\bcdot\})$ is finite
dimensional. Given any non-zero tripotent $e$, it follows from
finite-dimensionality and the repeated application of Result~\ref{R:po} that $e$ 
can be written as a sum of mutually orthogonal primitive tripotents. This brings us
to the final concept in this primer: the \emph{rank of a tripotent $e$} is the
minimum number of primitive tripotents required for such a decomposition of $e$ while the
\emph{rank of $(V,\{\bcdot,\bcdot,\bcdot\})$} is the highest rank that a tripotent of $V$ can have.
\medskip

\section{The boundary geometry of bounded symmetric domains}\label{S:bdry}

In this section we describe the boundary of a bounded symmetric domain in terms of the positive
Hermitian Jordan triple system associated to it. Thus, we shall
follow the notation introduced in Section~\ref{S:primer}. Recall that a bounded symmetric
domain $D$ has a realization as a bounded convex balanced domain. When we say ``Hermitian
Jordan triple system associated to $D$'', it is implicit that $D$ is this realization
and the association is the one given by \eqref{E:corrJTS}. This is a one-to-one
correspondence, described as follows:

\begin{result}[\cite{loos1977:bsdJp}, Theorem~4.1]\label{R:corresp}
 Let $D$ be a realization of a bounded symmetric domain as a bounded convex balanced domain
 in $\Cn$ for some $n\in \Z_+$. Then, $D$ is the open unit ball in $\Cn$ with respect
 to the spectral norm determined by the PHJTS associated to $D$. Conversely, given a
 PHJTS $(\Cn,\{\bcdot,\bcdot,\bcdot\})$, the open unit ball with respect to the spectral
 norm determined by it is a bounded symmetric domain $D$, and the PHJTS associated to $D$
 by the rule \eqref{E:corrJTS} is $(\Cn,\{\bcdot,\bcdot,\bcdot\})$.
\end{result}

\noindent In what follows, whenever we mention a bounded symmetric domain
$D$, {\em it will be understood that $D$ is a bounded convex balanced realization}.
\smallskip

The boundary of a bounded symmetric domain $D\subset \Cn$ has a certain stratification into
real-analytic submanifolds that can be described in terms of the PHJTS associated to $D$. The
first part of this section is devoted to describing this stratification. Fix a bounded symmetric
domain $D\subset \Cn$ and let $(\Cn,\{\bcdot,\bcdot,\bcdot\}_D)$ be the PHJTS associated to it.
It turns out (see \cite[Theorem~5.6]{loos1977:bsdJp}) that the set $M_D$ of tripotents
of $\Cn$ with respect to $\{\bcdot,\bcdot,\bcdot\}_D$ is a disjoint union of real-analytic
submanifolds of $\Cn$. For each $e\in M_D$, let $M_{D,e}$ denote the connected component of
$M_D$ containing $e$. The tangent
space $T_e(M_{D,e})$, viewed {\em extrinsically} (i.e., so that $e+T_e(M_{D,e})$ is the affine subspace of
all tangents to $M_{D,e}$ at $e$), is:
\[
 T_e(M_{D,e}) =  iA(e)\oplus V_1(e),
\]
where $A(e)$ is determined by the relation $V_2(e) = \{x+iy\in \Cn: x,y \in A(e)\}$, and
$V_j(e)$ is the eigenspace of $j=0,1,2$ in the Pierce decomposition of $\Cn$ with respect
to $e$. 
\smallskip

Let $M^*_D$ be the set of all non-zero tripotents and let $\|\bcdot\|_D$ denote the spectral
norm determined by $\{\bcdot,\bcdot,\bcdot\}_D$. Define
\begin{align*}
 E_D &:= \{(e,v)\in \Cn\times\Cn: e\in M^*_D \ \text{and} \ v\in V_0(e)\}, \\
 \bund_D &:= \{(e,v)\in E_D : \|v\|_D < 1\}.
\end{align*}
We can write $\bund_D$ as a disjoint union of the form
\begin{equation}\label{E:stratB}
 \bund_D \ := \ \bigsqcup_{\alpha\in \smoo}\bund_{D,\alpha},
\end{equation}
where $\smoo$ is the set of connected components of $M^*_D$, and each 
$\bund_{D,\alpha}$ is a connected, real-analytic submanifold of $\Cn\times \Cn$
that is a real-analytic fibre bundle whose fibres are unit $\|\bcdot\|_D$-discs. The key theorem about the
boundary of $D$ is as follows:

\begin{result}[\cite{loos1977:bsdJp}, Chapter~6]\label{R:strat}
 Let $D$ be a bounded symmetric domain in $\Cn$ and let $\eff: 
 \bund_D\to \Cn$ be defined by $\eff(e,v) := e+v$. Then:
 \begin{itemize}
  \item[$i)$] $\eff|_{\bund_{D,\alpha}}$ is an imbedding for each $\alpha\in \smoo$;
  \item[$ii)$] $\bdy D = \sqcup_{\alpha\in \smoo}\str_{D,\alpha}$, where
  $\str_{D,\alpha} := \eff(\bund_{D,\alpha})$;
  \item[$iii)$] in the above stratification of $\bdy D$, if $\str_{D,\alpha}$ is of
  dimension $d_\alpha$, then it is a closed, connected, real-analytic imbedded submanifold
  of the open set
  \[
   \Cn\setminus\!\!\!\!\bigcup_{\beta\,:\,\dim_{\R}(\str_{D,\beta}) < d_{\alpha}}\!\!\!\!\str_{D,\beta}.
  \]
 \end{itemize}
\end{result}

Furthermore, when $D$ is an {\em irreducible} bounded symmetric domain in $\Cn$, then we can
provide further information. Here, the rank of a bounded symmetric domain is the rank of the
Jordan triple system $(\Cn,\{\bcdot,\bcdot,\bcdot\}_D)$. 

\begin{result}[\cite{loos1977:bsdJp}, Chapter~6; \cite{vigue1991:domainesSymetriques}, 
Th{\'e}or{\`e}me~7.3]\label{R:stratIrr}
 Let $D$ be an irreducible bounded symmetric domain in $\Cn$ of rank $r$, and let $\smoo$ denote
 the set of connected components of $\bund_D$. Then, we have the following:
 \begin{itemize}
  \item[$i)$] $\smoo$ has cardinality $r$.
  \item[$ii)$] Each connected component of the decomposition \eqref{E:stratB} is a
  bundle over a submanifold of non-zero tripotents of rank $j, \ j\in \{1,\dots,r\}$. Denoting
  this bundle as $\bund_{D,\,j}, \ j\in \{1,2,\dots,r\}$, 
  we can express the stratification of $\bdy D$ given by Result~\ref{R:strat}-$(ii)$ as
  \[
   \bdy D = \bigsqcup_{j=1}^r\str_{D,\,j},
  \]
  where $\str_{D,\,j} := \eff(\bund_{D,\,j})$, and each $\str_{D,\,j}$ is connected.
  \item[$iii)$] The stratum $\str_{D,1}$ is dense in $\bdy D$.  
 \end{itemize}
\end{result}

The other goal of this section is to describe the structure of the germs of complex-analytic
varieties contained in the boundary of a bounded symmetric domain $D$. This structure can be
described in extremely minute detail; see, for instance, \cite{wolf72:fineStructure} by Wolf.
In fact, the papers about higher-rank bounded symmetric domains mentioned in 
Section~\ref{S:intro} make extensive use of this fine structure. However, {\em in this work, 
we only need very coarse information about the complex analytic structure of} $\bdy D$; specifically:
the distinction between the Bergman--Shilov boundary of $D$ and its complement in $\bdy D$. 
\smallskip

We denote the Bergman--Shilov boundary of $D$ by $\bdy_S D$.
We shall not formally define here the notion of the Shilov boundary of a uniform algebra; we
shall merely state that the Bergman--Shilov boundary of a bounded domain $D\Subset \Cn$ is 
the Shilov boundary of the uniform algebra $A(D) := \hol(D)\cap \smoo(\overline{D})$. However,
we do carefully state the following definition:

\begin{definition}
 Let $D$ be a bounded domain in $\Cn$.
 An {\em affine $\bdy D$-component} is an equivalence class under the equivalence relation
 $\thicksim_A$ on $\bdy D$ given by
 \[
  x\,\thicksim_A\,y \ \iff \ \text{$x$ and $y$ can be joined by a chain of segments lying in $\bdy D$},
 \]
 where a segment is a subset of $\Cn$ of the form $\{u+tv : t\in (0,1)\}, \ u, v\in \Cn$.
 A {\em holomorphic arc component of $\bdy D$} is an equivalence class under the equivalence relation
 $\thicksim_H$ on $\bdy D$ given by
 \[
  x\,\thicksim_H\,y \ \iff \ \text{$x$ and $y$ can be joined by a chain of analytic discs lying in $\bdy D$}.
 \]
\end{definition}

\noindent Roughly speaking, given a bounded domain $D\Subset \Cn$ and a point $x\in \bdy D$,
the holomorphic arc component of $\bdy D$ containing $x$ is the largest (germ of a) complex-analytic
variety lying in $\bdy D$ that contains $x$. The information that we require about holomorphic boundary
components is:
\begin{result}[\cite{loos1977:bsdJp}, Theorem~6.3]\label{R:BdyCmp}
 Let $D$ be the realization of a bounded symmetric domain as a
 bounded convex balanced domain in $\Cn$.
 \begin{itemize}
  \item[$i)$] The affine $\bdy D$-components and the holomorphic arc components of $\bdy D$ coincide.
  \item[$ii)$] A boundary component containing a point $x\in \bdy D$ is a non-empty open region in
  some $\C$-affine subspace of positive dimension passing through $x$ {\em unless} $x$ is a maximal tripotent.
 \end{itemize}
\end{result}

Finally, we mention the following description of the Bergman--Shilov boundary of a bounded symmetric
domain:

\begin{result}[\cite{loos1977:bsdJp}, Theorem~6.5]\label{R:BdyShi}
 Let $D\Subset \Cn$ be as in Result~\ref{R:BdyCmp}.
 The Bergman--Shilov boundary of $D$ coincides with each of the following sets:
 \begin{itemize}
  \item[$i)$] the set of maximal tripotents of $\Cn$ with respect to $\{\bcdot,\bcdot,\bcdot\}_D$;
  \item[$ii)$] the set of extreme points of $\overline{D}$;
  \item[$iii)$] the set of points of $\overline{D}$ having the maximum Euclidean distance from $0\in \Cn$.
 \end{itemize}
\end{result}
\medskip

\section{Some essential propositions}\label{S:essn}

This section contains several lemmas and propositions --- some being simple consequences of known
results, and some requiring substantial work --- that will be needed to prove our theorems. We begin
with the following result of Bell:

\begin{result}[\cite{bell1982:proper-circular}, Theorem~2]\label{R:bellsExt}
 Suppose $f: D_1\to D_2$ is a proper holomorphic map between bounded circular domains.
 Suppose further that $D_2$ contains the origin and that the Bergman kernel $K(w,z)$ associated to
 $D_1$ is such that for each compact subset $G$ of $D_1$, there is an open set $U = U(G)$ containing
 $\overline{D}_1$ such that $K(\bcdot, z)$ extends to be holomorphic on $U$ for each $z\in G$. Then
 $f$ extends holomorphically to a neighbourhood of $\overline{D}_1$.
\end{result}

Now let $D$ be any bounded balanced domain (not necessarily convex) in $\Cn$. If $D$ is not convex,
it will not be a unit ball with respect to some norm on $\Cn$. However, we do have a function that has
the same homogeneity property as a norm, with respect to which $D$ is the ``unit ball''. The
function $M_D: \Cn\to [0,\infty)$ defined by
\[
 M_D(z) := \inf\{t > 0 : z/t\in D\}
\]
is called the {\em Minkowski functional for $D$}. Assume that the intersection of each complex line
passing through $0\in \Cn$ with $\bdy{D}$ is a circle. Let $G$ be a compact subset of $D$. Then, as
$M_D$ is upper semicontinuous, $\exists r_G\in (0,1)$
such that $G\subset \{z\in \Cn : M_D(z) < r_G\}$
and the latter is an open set. Hence $z/r_G\in D \ \forall z\in G$. Clearly, $r_Gw\in D \ 
\forall w\in \{z\in \Cn : M_D(z) < 1/r_G\} =: U(G)$. By our assumptions, $\overline{D}\subset U(G)$.
 Let $K_D$ be the Bergman kernel of $D$. We recall that:
\[
 K_D(w,z) = \sum_{\nu\in \N}\psi_\nu(w)\overline{\psi_\nu(z)} \quad \forall (w,z)\in D\times D,
\]
where the right-hand side converges absolutely and uniformly on any compact subset of
$D\times D$ and $\{\psi_\nu\}_{\nu\in \N}$ is {\em any} complete orthonormal system for
the Bergman space of $D$ . Then --- owing to the fact that the
collection $\{C_\alpha z^\alpha : \alpha\in \N^n\}$ (where $C_\alpha > 0$ are suitable
normalization constants) is a complete orthonormal system for the Bergman space of $D$ ---
we can infer two things. First: the functions
\begin{equation}\label{E:extns}
 \phi_z(w) := K_D(r_Gw, z/r_G), \; \; w\in U(G),
\end{equation}
are well-defined by power series for each $z\in G$. Secondly:
\[
 K_D(r_Gw, z/r_G) = K_D(w,z) \; \; \forall (w,z)\in D\times G.
\]
Comparing this with \eqref{E:extns}, we see that each $\phi_z$ extends $K_D(\bcdot, z)$
holomorphically. In view of Result~\ref{R:bellsExt}, we have just deduced:

\begin{lemma}\label{L:ExtBell}
 Let $f: D_1\to D_2$ be a proper holomorphic map between bounded circular domains. Suppose
 $D_1$ and $D_2$ are both balanced. Assume that the intersection of every complex line
 passing through $0$ with $\bdy{D}_1$ is a circle. Then $f$ extends holomorphically to a
 neighbourhood of $\overline{D}_1$.
\end{lemma}

\noindent We remark that the above conclusion also follows from a later work
\cite{bell1993:algebraic-circular} of Bell.
\smallskip

Let $D$ be a bounded symmetric domain in $\Cn$. Let $(\Cn, \{\bcdot,\bcdot,\bcdot\}_D)$
be the Jordan triple system associated to $D$ (as in other places in this paper, we assume that
$D$ is a Harish-Chandra realization). Let $\dee_D$ and $Q_D$ be the maps \eqref{E:ops}
for the triple product $\{\bcdot,\bcdot,\bcdot\}_D$. We define the linear operators
$\ber{D}(x,y): \Cn\to \Cn$:
\[
 \ber{D}(x,y) := {\sf id}_D - \dee_D(x,y) + Q(x)Q(y), \; \; x, y\in \Cn.
\]
Consider the sesquilinear form $(x,y)\longmapsto {\rm Tr}[\dee_D(x,y)]$ on $\Cn$.
It turns out that the positivity of $\{\bcdot,\bcdot,\bcdot\}_D$ is equivalent to the above
sesquilinear form being an inner product on $\Cn$; see \cite[Chapter 3]{loos1977:bsdJp}. Furthermore
with respect to this inner product, we have:
\[
 \ber{D}(x,y)^* =\,\ber{D}(y,x) \; \; \forall x, y \in \Cn.
\] 
It is now easy to deduce that $\ber{D}(a,a)$ is a self-adjoint, positive semi-definite linear operator.
Consequently, $\ber{D}(a,a)$ admits a unique positive semi-definite square root, which we denote by
$\ber{D}(a,a)^{1/2}$. Having made these two definitions, we can state the following useful facts
about the geometry of $D$.

\begin{result}[\cite{loos1977:bsdJp}, Proposition~9.8; \cite{roos2000:JTS},
  Proposition~III.4.1]\label{R:AutFor}
 Let $D$ be the realization of a bounded symmetric domain as a convex balanced domain
 in $\Cn$. Fix a point $a\in D$ and let
 \[
  g_a(z) := a + \ber{D}(a,a)^{1/2}({\sf id}_D + \dee_D(z,a))^{-1}(z) \; \; \forall z\in D.
 \]
 Then, $g_a\in {\rm Aut}(D)$, $g_a(0) = a$, and 
 $g^\prime_a(z) = \ber{D}(a,a)^{1/2}\circ\ber{D}(z, -a)^{-1}$. Furthermore, $g_a^{-1} = g_{-a}$.
\end{result} 

Various versions of the following lemma have been known for a long time. We refer the reader
to \cite[Lemma\,15.2.2]{rudin1980:FT} for a proof.

\begin{lemma}\label{L:conv-const}
 Let $D$ be a bounded domain in $\Cn$ and let $p\in \bdy D$. Assume that there exists
 a ball $B$ centered at $p$ and a function $h\in \hol(B\cap D)\cap \smoo(\overline{B\cap D}; \C)$
 such that $h(p)=1$ and $|h(z)| < 1 \ \forall z\in \overline{B\cap D}\setminus\{p\}$.
 Let $a_0\in D$ and $\{\phi_k\}$ be a sequence of automorphisms of
 $D$ such that $\phi_{k}(a_0)\lrarw p$ as $k\to \infty$. Then, $\{\phi_k\}$ converges
 uniformly on compact subsets of $D$ to $\con{p}$ --- the map that takes the constant value $p$.
\end{lemma}

We now state a version of Schwarz's lemma for convex balanced domains and
then a version of Schwarz's lemma for irreducible bounded symmetric
domains, both of which are needed in the proof of our Key Lemma (see Section~\ref{S:intro}).

\begin{result}[\cite{rudin1980:FT}, Theorem~8.1.2] \label{R:SL-cbd}
  Let $\Omega_1$ and $\Omega_2$ be balanced regions in $\Cn$
  and $\C^m$ respectively, and let $F: \Omega_1\to \Omega_2$ be a holomorphic
  map. Suppose $\Omega_2$ is convex and bounded. Then:
 \begin{enumerate}
 \item[$i)$] $F ^{\prime}(0)$ maps $\Omega_1$ into $\Omega_2$; and
 \item[$ii)$] $F(r\Omega_1) \subseteq r\Omega_2 \ (0 < r \leq 1)$ if $F(0) = 0$.
 \end{enumerate}
\end{result}

\begin{result}[\cite{vigue1991:domainesSymetriques}, Th{\'e}or{\`e}me~7.4] \label{R:SL-vig}
  Let $D$ be an irreducible bounded symmetric domain in $\Cn$ in its Harish-Chandra
  realization (whence it is the unit ball in $\Cn$ for the associated spectral norm
  $\|\bcdot\|$). Let $F:D \to D$ be a holomorphic
  map such that $F(0) = 0$. Assume that for some non-empty open set $U \subset D$, we
  have $\|F(z)\| =\|z\| \ \forall z \in U$. Then $F$ is an automorphism of $D$.
\end{result}

With these two results, we can now give a proof of our Key Lemma:
\smallskip

\noindent{\em The proof of the Key Lemma.} Let $z \in U$, and set $w := F(z)$. By
hypothesis, $F$ maps $\Delta_z$ into $D$ and $\left(F|_{W_1}\right)^{-1}$ maps $\Delta_w$ into
$D$. Applying Result~\ref{R:SL-cbd} to $F|_{\Delta_z}$ and to
$\left(F|_{W_1}\right)^{-1}\big|_{\Delta_w}$, we have
$\|F(z)\| = \|z\|$ for every $z \in U$. Thus by the Schwarz lemma for irreducible bounded
symmetric domains, $F$ is an automorphism of $D$.  \qed
\smallskip

We now state and prove a technical proposition regarding the invertibility of the operator
$\ber{D}$ at certain off-diagonal points in $\bdy D\times \bdy D$, where $D$ is an {\em irreducible}
bounded symmetric domain of dimension\,$\geq 2$. Here $\mathcal{M}_{D,1}$ denotes the stratum of
$\bdy D$ described by Result~\ref{R:stratIrr}. This result and our Key Lemma are the central
ingredients in the proof of our Main Theorem.

\begin{proposition}\label{P:InvLem}
  Let $D$ be the realization of an irreducible bounded symmetric domain of dimension $n$ as a bounded
  convex balanced domain in $\Cn$, $n\geq 2$. Let $p \in \bdy D$. For each $z_0 \in \mathcal{M}_{D,1}$
  and each $\mathcal{M}_{D,1}$-open neighbourhood $U \ni z_0$, there exists a point $w \in U$ such
  that $\det\ber{D}(\bcdot,p)$ is non-zero on the set $\{\zeta w:\zeta \in \C, |\zeta| = 1\}$.
\end{proposition}
 
\begin{remark}\label{rem:maxStratum}
 In the following proof, we argue by assuming that the conclusion above is false. We can
 instantly arrive at a contradiction at the point $(\bullet)$ in the proof below if we invoke results
 on the fine structure of $\bdy D$; see \cite{koranyiWolf1965:symmDomains} or 
 \cite{wolf72:fineStructure}, for instance. However, we
 provide an elementary argument beyond $(\bullet)$ to complete the proof in the hope that
 appropriate analogues of the above may be formulated in other contexts.
\end{remark}  
\begin{proof}
Let us denote $\det\ber{D}(z,p)$ as $h(z)$, where $z \in \Cn$. Let us assume
that the result is false. Then, there exists a point $z_0 \in
\mathcal{M}_{D,1}$ and an $\mathcal{M}_{D,1}$-open neighbourhood $U \ni z_0$ such
that for each $w \in U$, there exists a $\zeta_w \in \{\zeta \in \C : |\zeta| = 1\}$ with
$h(\zeta_ww) = 0$. Let $q$ denote the quotient map $q:\Cn
\setminus \{0\} \to \pro^{n-1}$. Also write
\[
  Z_h := h^{-1}\{0\}, \quad Z := Z_h \cap \mathcal{M}_{D,1}.
\]
Our assumption implies that $q(Z)$ contains a non-empty open set $\mathcal{V}
\subset \pro^{n-1}$. Let $\mathcal{A} := \{z \in \Cn : 1 - \eps < \|z\| < 1 + \eps \}$,
where $\|\bcdot\|$ denotes the spectral norm relative to which $D$ is the unit ball, and
$\eps$ is a {\em fixed} number in $(0,1)$. As $\mathcal{V}\subset q(\mathcal{A})$, it is easy to
see that $\mathcal{V}$ can be covered by
finitely many holomorphic coordinate patches $(U_1,\psi_1),\dots,(U_M,\psi_M)$
such that the maps
\[
  q_j := \psi_j \circ q|_{q^{-1}(U_j) \cap \mathcal{A}} : q^{-1}(U_j) \cap
  \mathcal{A} \to \C^{n-1}
\]
are Lipschitz maps. Since Lipschitz maps cannot increase Hausdorff
dimension (see \cite[Proposition~14.4.4]{rudin1980:FT}, for instance)
and $\text{dim}_{\R}(\mathcal{V}) = 2n - 2$, the preceding
discussion shows that the Hausdorff dimension of $Z$ (and hence the dimension of
$Z$ as a real-analytic set) is $2n - 2$. As
$Z_h$ is a complex analytic subvariety, its singular locus is of complex
dimension $\leq n - 2$. Thus, we can find a point $x_0 \in Z$ that is a
regular point of $Z_h$, and an open ball $B$ around $x_0$ that is so small that
\begin{itemize}
 \item $\mathcal{M}_{D,1}\cap B$ is a submanifold of $B$;
 \item $B \cap Z_h$ is an $(n-1)$-dimensional complex submanifold of $B$;
 \item the dimension of $B \cap Z$ is $2n -2$.
\end{itemize}
These three facts imply that $M : = B \cap Z_h \subset \mathcal{M}_{D,1}$. We can 
deduce this by considering a local defining function $\rho_B : B\to \R$ for
$\mathcal{M}_{D,1}$ and observing that, by {\L}ojasiewicz's theorem
\cite{Lojasiewicz1959:probDiv}, $\left.\rho_B\right|_M \equiv 0$. {\em If $D = \mathbb{B}^n$,
we already have a contradiction} and, hence, the proof.
\smallskip

Since $\mathcal{M}_{D,1}$ is a real-analytic submanifold of
$\Cn \setminus \bigsqcup_{j \geq 2}\mathcal{M}_{D,\,j}$, where $\mathcal{M}_{D,\,j}$
are the strata of $\bdy D$ discussed in Section~\ref{S:bdry}, we can define the Levi-form of
$\mathcal{M}_{D,1}$ --- denoted by $\mathfrak{L}(z,V), \ z \in \mathcal{M}_{D,1}, \
V \in H_z(\mathcal{M}_{D,1})$. A few words about notation: in this proof, we shall 
work with the tangent bundle of $\mathcal{M}_{D,1}$ defined extrinsically. So, when
referring to vectors in $T_z(\mathcal{M}_{D,1})$, we shall view them either as real or
as complex vectors, as convenient, such that $z + T_z(\mathcal{M}_{D,1})$ is the
hyperplane tangent to $\mathcal{M}_{D,1}$ at $z \in \mathcal{M}_{D,1}$. In
this scheme:
 \[
   H_z(\mathcal{M}_{D,1}) := T_z(\mathcal{M}_{D,1}) \cap
  iT_z(\mathcal{M}_{D,1}).
\]
As $\text{dim}_{\C}(M) = n - 1, \ \mathfrak{L}(z,\bcdot) \equiv 0 \ \forall z \in
M$. The curve $\gamma : (-\eps,\eps) \to \mathcal{M}_{D,1}$ (for $\eps > 0$
suitably small) $\gamma(t) := \exp(it)z$ is transverse to $M$ at $z$. This is
because if $\gamma'(0) = iz$ were in $H_z(\mathcal{M}_{D,1})$, then 
\[
  i\gamma'(0) = -z \in H_z(\mathcal{M}_{D,1}) \subset T_z(\mathcal{M}_{D,1}),
\]
which contradicts the convexity of $D$. Consequently, for $\eps_0 > 0$
sufficiently small, the set $\{\exp(it)z : t \in (-\eps_0,\eps_0), z \in M \}$
contains an $\mathcal{M}_{D,1}$-open neighbourhood of $x_0$. Thus, $\mathcal{M}_{D,1}$
is Levi-flat at $x_0$. As $\mathcal{M}_{D,1}$ is real-analytic, it is a Levi-flat
hypersurface. 
\smallskip

We shall now show that Levi-flatness of $\mathcal{M}_{D,1}$ leads to a
contradiction. Let us pick an $x \in \mathcal{M}_{D,1}$. Owing to
Levi-flatness, we can find a ball $B_x$, centered at $x$, such that
\[
 D_x^- := D \cap B_x, \ \ \ D_x^+ := B_x \setminus \overline{D}
\]
are both pseudoconvex. Let $\boldsymbol{n}_x$ denote the unit outward normal
vector to $\bdy D$ at $x$ ($x \in \mathcal{M}_{D,1}$). Owing to convexity of $D$,
we can find an $\eps_0 > 0$ and a $\delta_0 > 0$ such that
\[
  H_x(\eps_0;\delta) := x +\delta\boldsymbol{n}_x + \{V \in
  H_x(\mathcal{M}_{D,1}) : |V| < \eps_0\} \subset D_x^+
\]
for each $\delta \in (0,\delta_0)$. Here, $|\bcdot|$ denotes the Euclidean
norm. As $H_x(\eps_0;\delta)$ is a copy of a complex $(n-1)$-dimensional ball
and as $D_x^+$ is taut --- see
\cite[Proposition~2.1]{kerzmanRosay1981:domainesTaut} --- it follows that 
$H_x(\eps_0;0) \subset \mathcal{M}_{D,1}$. To summarize, $\mathcal{M}_{D,1}$ has
the following property:
\begin{itemize}
 \item[$(\bullet)$] At each $x \in \mathcal{M}_{D,1}$, a germ of the set
 $(x+H_x(\mathcal{M}_{D,1}))$ lies in $\mathcal{M}_{D,1}$. 
\end{itemize}
\smallskip

Let us now pick and fix a point $y^0 \in \mathcal{M}_{D,1}$. Let
$(z_1,\dots,z_n)$ be global holomorphic coordinates in $\Cn$, associated to
an appropriate rigid motion of $D$, such that $y^0 = (0,\dots,0), D \subset
\{\re z_1 > 0  \}$ and $H_{y^0}(\mathcal{M}_{D,1}) = \{z_1 = 0\}$ relative to
these coordinates. Let $W$ be a non-zero vector in $H_{y^0}(\mathcal{M}_{D,1})$ and let 
$D_W := D \cap \text{span}_\C\{W,\boldsymbol{n}_{y^0}\}$. Clearly, $D_W$ is convex and by $(\bullet)$
$\mathcal{M}_{D,1} \cap \text{span}_{\C}\{W,\boldsymbol{n}_{y^0}\} =: \mathcal{M}_W$ has the
property that for each point $y \in \mathcal{M}_W$, the germ of a complex line through $y$, call
it $\Lambda_{y,W}$, lies in $\mathcal{M}_W$. Let us view $D_W$ as lying in $\C^2$,
whence a portion of $\mathcal{M}_W$ near $(0,0)$ can be parametrized by three real
variables as follows:
\[
  r(t,u,v) = \rho(t) + a(t)(u + iv), \; \; |t| < \eps_1, \ |u|, \ |v| < \eps_2,
\]
where $\rho:(-\eps_1,\eps_1) \to \mathcal{M}_W$ is a smooth curve through $(0,0)$ such
that $\rho'(t)$ is orthogonal to $\Lambda_{\rho(t),W}$ for each $t$, and
$a:(-\eps_1,\eps_1) \to \C^2$ is such that $a(t)$ is parallel to $\Lambda_{\rho(t),W}$
for each $t$. For the remainder of this paragraph, $\boldsymbol{n}(t,u,v)$ will denote
the {\em inward} unit normal to $\bdy D_W$ at $r(t,u,v)$, and $\bcdot$ will denote the
standard inner product on $\R^4$. Define the matrix-valued function $\Gamma : (-\eps_1,
\eps_1) \times (-\eps_2,\eps_2)^2  \to \R^{3 \times 3}$ by 
\[
  \Gamma(\tau,U,V) := \left.{\sf Hess}_{t,u,v}\left( r(t,u,v) \bcdot \boldsymbol{n}(\tau,U,V)
	\right)\right|_{(t,u,v) = (\tau,U,V)}\,.
\]
The convexity of $D_W$ implies that $\Gamma(\tau,U,V)$ is {\em positive} semidefinite at each
$(\tau,U,V)$ (recall that $\boldsymbol{n}(\tau,U,V)$ is the inward normal at $r(\tau,U,V)$).
By choosing $\eps_1,\eps_2 > 0$ small enough, we can ensure that
$(n_1^2 + n_2^2)(t,u,v) \neq 0$ for every $(t,u,v)$, where we
write $\boldsymbol{n} = (n_1,n_2,n_3,n_4)$, and that $a$ is of the form
$a(t) = (\alpha(t) + i\beta(t),1)$. We compute to observe that two of the principal
minors of $\Gamma$ turn out to be $-(n_1\alpha' + n_2\beta')^2$ and
$-(n_2\alpha' - n_1\beta')^2$, which must be non-negative. This gives us the system of
equations
\begin{align*}
  \left.n_1\alpha' + n_2\beta'\right|_{(\tau,U,V)}\,&=\,0 \\
  \left.-n_1\beta' + n_2\alpha'\right|_{(\tau,U,V)}\,&=\,0 \; \; \forall (\tau,U,V).
\end{align*}
By our assumption on $\boldsymbol{n}$, this implies that $\alpha' = \beta' \equiv 0$.
Restating this geometrically, there is a small $\mathcal{M}_W$-open neighbourhood of
$0 \in \bdy D_W$ such that, for every $y$ in this neighbourhood,
$\Lambda_{y,W}$ is parallel to the vector $W$. This holds true for each non-zero
$W \in H_{y^0}(\mathcal{M}_{D,1})$. Thus, there is an
$\mathcal{M}_{D,1}$-open patch $\omega \ni y^0$ such that 
\begin{equation}\label{E:fol}
  x + H_x( \mathcal{M}_{D,1}) \text{ is parallel to } \{z_1 = 0 \} \ \text{for every $x \in \omega$}.
\end{equation}

By Result~\ref{R:stratIrr}, $\mathcal{M}_{D,1}$ is connected. Thus, if $y^0\neq y\in \mathcal{M}_{D,1}$, 
then $y$ can be joined to $y^0$ by a chain of $\mathcal{M}_{D,1}$-open patches 
$\ome_0,\dots,\ome_N$, where $\ome_0$ equals the patch $\ome$ in \eqref{E:fol}, 
$\ome_{j-1}j\cap \ome_j\neq \varnothing, \ j = 1,\dots,N$, and $\ome_N\ni y$. By a standard
argument of real-analytic continuation, we deduce that \eqref{E:fol} holds with $\ome_N$ replacing
$\ome$ (where $z_1$ comes from the global system of
coordinates fixed at the beginning of the previous paragraph). Hence,
$x + H_x(\mathcal{M}_{D,1})$ is parallel to $\{z_1 = 0\}$ for each $x \in \mathcal{M}_{D,1}$.
As $\mathcal{M}_{D,1}$ is dense in $\bdy D$, and $D$ is bounded, we can find a $\xi \in D$ and a
vector $W = (W_1,\dots,W_n)$ with $W_1 = 0$ such that the ray $\{\xi + tW ; t \geq 0\}$
intersects $\bdy D$ at a point in $\mathcal{M}_{D,1}$. Then, this ray must be tangential to
$\mathcal{M}_{D,1}$ at the point of intersection, which is absurd as $D$ is convex.
Hence, our initial assumption must be false.
\end{proof}
\medskip

\section{The proof of Theorem~\ref{T:mThm}}\label{S:mThm}
Before we proceed further, we clarify our notation for the different norms that
will be used in the proof of Theorem~\ref{T:mThm}. With $D_1$ and $D_2$ as in
Theorem~\ref{T:mThm}, $\|\bcdot\|_j$ will denote the spectral norms such that $D_j$
is the unit $\|\bcdot\|_j$-ball in $\Cn, j =1,2$. The Euclidean norm on $\Cn$ will be
denoted by $|\bcdot|$. We will also need to impose norms on certain linear operators
on $\Cn$. We shall use the operator norm induced by the Euclidean norm: for a
$\C$-linear operator $A$ on $\Cn$, we set
\[
  \|A\|_{op} := \sup_{|x| = 1}|Ax|.
\]

\noindent{\bf The proof of Theorem~\ref{T:mThm}.} We shall take $D_1$ and $D_2$ to be
Harish-Chandra realizations of the given bounded symmetric domains. We may
assume, composing $F$ with suitable automorphisms if necessary, that $F(0) = 0$.
\smallskip

By Lemma~\ref{L:ExtBell}, $F$ extends to a holomorphic map defined on a
neighbourhood $N$ of $\overline{D}_1$. For simplicity of notation, we shall
denote this extension also as $F$. The complex Jacobian $\cjac F$ is holomorphic
on $N$ and $\cjac F \not\equiv 0$ on $D_1$. Hence, by the maximum principle, $\cjac
F \not\equiv 0$ on $\bdy D_1$. By definition, we can find a point $p \in \bdy_S
D_1$ such that 
\[
  \sup_{\overline{D}_1}|\cjac F| = |\cjac F(p)| \ \neq 0.
\]
By the inverse function theorem, we can find a ball $B(p,r) \subset N$ such
that $F|_{B(p,r)}$ is injective. Let us write
\[
  \OM_1 := B(p,r) \cap D_1, \; \; \OM_2 := F(B(p,r)) \cap D_2.
\]
We shall use our Key Lemma (see Section~\ref{S:intro}, and Section~\ref{S:essn}
for its proof) to deduce the result. The regions $W_1$ and $W_2$ of that Lemma
will be constructed by applying suitable automorphisms to $\OM_1$ and $\OM_2$.
\smallskip

\noindent{\bf Claim}. $F(p) \in \bdy_S D_2$.

\noindent Suppose $F(p) \not\in \bdy_S D_2$. It follows from Result~\ref{R:BdyCmp} and
Result~\ref{R:BdyShi} that there is a vector $V \in \Cn \setminus \{0\}$ and
neighbourhood $\ome$ of $0 \in \C$ such that $\psi(\ome) \subset F(B(p,r)) \cap \bdy
D_2$, where $\psi : \ome \ni \zeta \longmapsto F(p) + \zeta V$. Next, define
\[
  \widetilde{\psi} := (F|_{B(p,r)})^{-1} \circ \psi. 
\]
Since $F|_{D_1}$ is proper and $F|_{B(p,r)}$ is injective,
\[
  F(z) \in F(B(p,r)) \cap \bdy D_2 \iff z \in B(p,r) \cap \bdy D_1.
\]
Thus $\widetilde{\psi}(\ome) \subset \bdy D_1$. Furthermore,
$\widetilde{\psi}$ is non-constant and $\widetilde{\psi}(0) = p$. By definition,
each point of $\widetilde{\psi}(\ome) \setminus \{p\}$ lies in the holomorphic
arc component of $\bdy D_1$ containing $p$. This is a contradiction
since $p$, being an extreme point, is a one-point affine $\bdy D_1$-component and
thus, by Result~\ref{R:BdyCmp}, a one-point holomorphic arc component of $\bdy
D_1$. Hence the claim.
\smallskip

Let us now take a sequence $\{a_k\} \subset \OM_1$ such that $a_k \to p$, and
let $b_k := F(a_k)$. Let $\phi_k^1 \in \text{Aut}(D_1)$ denote an automorphism that
maps $0$ to $a_k$. Let $\phi_k^2 \in \text{Aut}(D_2)$ be an automorphism that maps $0$
to $b_k$. Owing to Result~\ref{R:BdyShi} and to convexity, we can construct a peak function
for $p$ on $\overline{D}_1$. Likewise (in view of the last claim) $F(p)$ is
a peak point of $D_2$. By Lemma~\ref{L:conv-const}, we get:   
\begin{equation}
  \phi_k^j \lrarw \con{p^j} \ \text{uniformly on compacts}, \; \; j =1,2,
  \label{E:CvgCon}
\end{equation}
where $p^j, p^1 : = p, p^2 := F(p)$.
\smallskip

We now define
\[
  \OM^k_j := (\phi_k^j)^{-1}(\OM_j), \; \; j =1,2, \ k \in \Z_+.
\]
Given any $r > 0$, write $rD_j := \{ z \in \Cn : \|z\|_j < r\}, \ j=1,2$.
By \eqref{E:CvgCon}, there exists a sequence $k_1 < k_2 < k_3 < \dots$ in $\Z_+$
such that 
\[
  \phi_{k_\nu}^1\left( (1 - 1/s)\overline{D}_1\right) \subset \OM_1 \; \; \forall
  \nu \geq s, \ s \in \Z_+.
\]
By \eqref{E:CvgCon} again, we can extract a sequence of indices $\nu(1) < \nu(2) <
\nu(3) < \dots$ such that
\[
  \phi_{k_{\nu(t)}}^2\left( (1 - 1/s)\overline{D}_2\right) \subset \OM_2 \; \; \forall
  t \geq s, \ s \in \Z_+.
\]
In the interests of readability of notation, let us re-index $\{k_{\nu(s)}\}_{s
\in \Z_+}$ as $\{k_m\}_{m \in \Z_+}$. Then, the above can be summarized as:
\begin{itemize}
  \item[$(*)$] With the sequences of maps $\{\phi_k^1\} \subset \text{Aut}(D_1)$ and
	$\{\phi_k^2\} \subset \text{Aut}(D_2)$ as described above, there is a sequence
	$\{k_m\}_{m \in \Z_+} \subset \Z_+$ and a strictly increasing
	$\Z_+$-valued function $\nu^*$ such that
	\begin{align}
	  (1 - 1/s)\overline{D}_1 &\subset \OM_1^{k_m} \;\; \forall m \geq s, \ s \in \Z_+,  \\
	  (1 - 1/\nu^*(s))\overline{D}_2 &\subset \OM_2^{k_m} \; \; \forall m \geq s, \ s \in \Z_+.
	  \label{E:Lsub}
	 \end{align}
\end{itemize}
\smallskip

\noindent{\bf Step 1.} \emph{Analysing the family $\{(\phi_{k_m}^2)^{-1} \circ
F \circ \phi_{k_m}^1\}_{m \in \Z_+}$}

\noindent Consider the maps $G_m : D_1 \to D_2$ defined by 
\[
  G_m := (\phi_{k_m}^2)^{-1} \circ F \circ \phi_{k_m}^1.
\]
By Montel's theorem, and passing to a subsequence and relabelling if necessary, we get a map
$G \in \hol(D_1;\Cn)$ such that $G_m \to G$ uniformly on compact subsets. Let us
fix an $s \in \Z_+$. By $(*)$, we infer that $\exists M_s \in \Z_+$ such that
$(1-1/s)\overline{D}_j \subset \OM_j^{k_m} \ \forall m \geq M_s$, $j = 1,2$.
Note that $G_m|_{\OM_1^{k_m}}$ is a biholomorphism, whence $G_m'(0)$ is
invertible for each $m$. Hence, by the
Schwarz lemma for convex balanced domains (i.e. Result~\ref{R:SL-cbd} above)
$G_m'(0)$ maps $(1-1/s)D_1$ into
$D_2$ and $G_m'(0)^{-1}$ maps $(1-1/s)D_2$ into $D_1 \ \forall m \geq M_s$. We
claim that this implies that $G'(0)$ is invertible. Suppose not. Then we would
find a $z_0$ with $\|z_0\|_1 = (1-2/s)$ such that $G'(0)z_0 = 0$. Note that
$G_m'(0) \to G'(0)$ in norm, whence, given any $\eps > 0, \|G_m'(0)z_0\|_2 <
\eps$ for every sufficiently large $m$. If we now choose $\eps \leq (1-2/s)^2$,
we see that
\[
  G_m'(0)^{-1}\left(\{ \|w\|_2 = (1-2/s)\}\right) \not\subset D_1
\]
for all sufficiently large $m$. This is a contradiction. Hence the claim.
\smallskip

Now that it is established that $G'(0)$ is invertible, it follows that
$G_m'(0)^{-1} \to G'(0)^{-1}$ in norm. Hence, $G'(0)^{-1}$ maps $(1-1/s)D_2$ into
$D_1$. Recall that $s \in \Z_+$ was arbitrarily chosen and that the function $\nu^*$
in $(*)$ is strictly increasing. Thus, $G'(0)^{-1}$ maps $D_2$ into $D_1$. By 
construction, $G(D_1) \subset \overline{D}_2$. Now, $D_2$ is complete (Kobayashi)
hyperbolic. Hence $D_2$ is taut; see \cite{kiernan1970:relations}.
As $G(0) = 0 \in D_2, G$ maps $D_1$ to $D_2$. So, the
holomorphic map $G'(0)^{-1} \circ G : D_1 \to D_1$ satisfies all the conditions
of Cartan's uniqueness theorem. Thus,
\[
  G'(0)^{-1} \circ G = {\sf id}_{D_1},
\]
which means that $G = G'(0)|_{D_1}$.
\smallskip

\noindent{\bf Step 2.} \emph{Showing that $D_1$ and $D_2$ are biholomorphically equivalent}

\noindent We have shown in Step 1 that $G'(0)^{-1}$ maps $(1-1/s)D_2$ into $D_1$. As
$G'(0)$ is injective, this means that $G'(0)(D_1)$ contains $(1-1/s)D_2$ for arbitrarily large
$s\in \Z_+$. Thus $G$ maps $D_1$ onto $D_2$. It follows that $D_1$ is biholomorphic to
$D_2$.
\smallskip

It would help to simplify our notation somewhat. By the
nature of the argument in Step 1, it is clear that we can assume that the
sequences $\{a_k\} \subset \OM_1$ and $\{b_k\} \subset \OM_2$ are so selected
that $(*)$ is true with $\{k_m\}_{m \in \Z_+} = \{1,2,3,\dots\}$. Owing to Step 2, we
may now assume $D_1 = D_2 := D$.  The argument we will make in Step 3 below
is valid regardless of the specific sequence $\{a_k\}$ or $\{b_k\}$. Hence, in the next
three paragraphs following this, the sequence $\{A_k\}$ will stand for either
$\{a_k\}$ or $\{b_k\}$, and the point $q$ will stand for either $p$ or $F(p)$.
Also, we will abbreviate $\phi^{j}_{A_k}$ to $\phi_k$. 
\medskip

\noindent{\bf Step 3.} \emph{Producing subsequences of $\{\phi_k\}$
that converge on ``large'' subsets of $\bdy D$}.

\noindent By Result~\ref{R:AutFor} we may take $\phi_k = g_{A_k}$, whence
\begin{equation}\label{E:phiJac}
  \phi_k'(z) = \ber{D}(A_k,A_k)^{1/2} \circ \ber{D}(-z,A_k)^{-1}.
\end{equation}
In the argument that follows, it is implicit that each $\phi_k$ is defined as a
holomorphic map on some neighbourhood (which depends on $\phi_k$) of
$\overline{D}$; see Lemma~\ref{L:ExtBell}. By Proposition~\ref{P:InvLem}
we can find a point $\xi_0 \in \mathcal{M}_{D,1}$ such that
\[
  \det\ber{D}(e^{i\theta}\xi_0,q) \neq 0 \; \; \forall \theta \in \R.
\]
By continuity, there exists a $\overline{D}$-open neighbourhood $\Gamma$ of $q$,
an $\mathcal{M}_{D,1}$-open neighbourhood $W$ of $\xi_0$, and a
$\overline{D}$-open set $V$ with the following properties:
\begin{itemize}
  \item[(a)] $z \in V \implies e^{i\theta}z \in V \ \forall \theta \in \R$;
  \item[(b)] $V \cap \bdy D = S^1 \bcdot W$;
  \item[(c)] $z \in V \implies tz \in V \  \forall t \in [1,1/\|z\|]$
\end{itemize}
(now $\|\bcdot\|$ is the spectral norm associated to $D$); such that
\begin{equation}\label{E:BerInv}
  \det \ber{D}(z,w) \neq 0 \; \; \forall (z,w) \in \overline{V} \times \Gamma.
\end{equation}
Here, given a set $X \subset \Cn, S^1 \bcdot X$ stands for the set
$\{e^{i\theta}x : x \in X, \theta \in \R\}$. Let us call any pair $(V,W)$, where
$V$ is a $\overline{D}$-open set and $W$ is an $\mathcal{M}_{D,1}$-open set, a
{\em truncated prism with base $S^1 \bcdot W$} if $(V,W)$ satisfies properties
(a)-(c) above.

We can find $V'$ and $W'$, with $\overline{W'} \subset W$, such that $(V',W')$
is a truncated prism with base $S^1 \bcdot W'$ with the properties:
\begin{itemize}
  \item $\overline{V'} \subset V$;
  \item There exists a $\delta_0 \ll 1$ such that for $z_1,z_2 \in
	\overline{V'}$, the segment $[z_1,z_2] \subset V$ whenever $|z_1 - z_2| <
	\delta_0$.
\end{itemize}
Owing to holomorphicity and convexity,
\begin{equation}\label{E:MeanVal}
   \phi_k(z_1) - \phi_k(z_2) = \int_0^1 \phi_k' \left(z_1 + t(z_2 -
  z_1)\right)(z_2 - z_1) dt, \; \; z_1,z_2 \in \overline{D}.
\end{equation}
We can find a $K \equiv K(W)$ such that, in view of \eqref{E:BerInv},
$\{\ber{D}(z,A_k): k \geq K(W), z \in \overline{V}\}$ is a compact family in
$GL(n,\C)$. Hence, in view of \eqref{E:phiJac} (and since 
$\{\ber{D}(A_k,A_k): k\in \Z_+\}$ is a relatively compact family in $\C^{n\times n}$),
there exists a constant $C > 0$ such
that
\[
  \|\phi_k'(z)\|_{op}\,\leq\,C \ \forall z \in \overline{V}, \; \; \forall k \geq K.
\]
By our construction of $V'$, and from \eqref{E:MeanVal}, we conclude:
\[
  |\phi_k(z_1) - \phi_k(z_2)| \leq C|z_1 - z_2| \; \; \forall z_1,z_2 \in
  \overline{V'}, \ |z_1 - z_2| < \delta_0, \ \text{and} \ \forall k \geq K.
\]
In short, $\{\phi_k|_{\overline{V'}}\} \subset \mathcal{C}(\overline{V'};\Cn)$ is
an equicontinuous family.
\smallskip

By the Arzela-Ascoli theorem, we can find a subsequence of $\{\phi_k\}$ that
converges uniformly to $q$ on $\overline{V'}$. For simplicity of notation, let
us continue to denote this subsequence as $\{\phi_k\}$. Then there exists a $K_1
\in \Z_+$ such that $\phi_k(\overline{V'}) \subset \OM$ (which denotes either
$\OM_1$ or $\OM_2$) $\forall k \geq K_1$. Furthermore, we may assume that $K_1$
is so large that, thanks to $(*)$,
\[
  (1-1/s)\overline{D} \subset \phi_k^{-1}(\OM) \; \; \forall k \geq K_1,
\]
where $s$ is so large that $(1-1/s)\overline{D} \cap V'$ is a
non-empty open set. By construction:
\[
  z \in V' \cap D \implies \Delta_z \subset (1-1/s)\overline{D} \cup V'.
\]
Hence $\Delta_z \subset \phi_k^{-1}(\OM) \ \forall k \geq K_1$. We summarize
the content of this paragraph as follows:
\begin{itemize}
  \item[$(**)$] Given any truncated prism $(V,W)$ with base $S^1 \bcdot W$ such
	that $\ber{D}(z,A_k) \neq 0$ on $\overline{V}$ for all $k$ sufficiently
	large, we can find a $K_1 \in \Z_+$ and a truncated prism $(V',W')$ with
	$\overline{V'} \subset V$ such that $\Delta_z \subset \phi_k^{-1}(\OM)$ for each
	$z \in V' \cap D$ and each $k \geq K_1$.
\end{itemize}
\smallskip

\noindent{\bf Step 4.} \emph{Completing the proof.}

\noindent By Proposition~\ref{P:InvLem} and $(**)$, we can find a truncated
prism $(V',W')$ with base $S^1 \bcdot W'$ which has all the properties stated in
$(**)$. Let $s \in \Z_+$ be so large that $(1-1/s)D \cap V' := U'$
is a non-empty open set. As $G_k \to G$ uniformly on $U'$ (by Step 1), there exists a point
$w_0 \in G(U')$, $K_2 \in\Z_+$ and a $c > 0$ such that the ball
\[
  B(w_0,c) \subset G(U') \cap G_k(U') \; \; \text{and} \; \; B(w_0,c) \subset \OM_2^k \; \;
  \forall k \geq K_2.
\]
Write $\|\bcdot\|$ for the spectral norm associated to $D$. Let
$R: \Cn \setminus {0} \to \bdy D$ be given by $R(w) := w/\|w\|$. By
Proposition~\ref{P:InvLem} and $(**)$, we can find a $\mathcal{M}_{D,1}$-open
subset $\ome_2$ such that
\[
  \ome_2 \subset R(B(w_0,c)),
\]
a truncated prism $(V_2,\ome_2)$ with base $S^1 \bcdot \ome_2$, and a $K_3 \in
\Z_+$ such that $\Delta_w \subset \OM_2^k$ for each $w \in V_2 \cap D$ and
each $k \geq K_3$. Let us now set $U := G^{-1}\left(R^{-1}(\ome_2) \cap B(w_0,c)\right)$,
and $K^* := \max(K_1,K_2,K_3)$. Finally, we set
\[
  W_j := \left(\phi^j_{K^*}\right)^{-1}\!\!(\OM_j), \; \; j = 1,2,
\]
with the understanding that $\phi^1_k = g_{a_k}$ and $\phi^2_k = g_{b_k}$.
\smallskip

As $U\subset V', \Delta_z \subset W_1$ for each $z \in U$. By construction 
\[
  G_{K^*}(z) \in B(w_0,c) \subset W_2 \; \; \forall z \in U.
\]
Finally, by construction, for each $z\in U$, there exists a point $w_z \in
\Delta_{G_{K^*}(z)}$ that belongs to $V_2 \cap D$. Thus, $\Delta_{G_{K^*}(z)}
\subset W_2$. Recall that $\left.G_{K^*}\right|_{W_1} : W_1 \to W_2$ is a biholomorphism and $G_{K^*}(0)
= 0$. By our Key Lemma, $G_{K^*}$, and consequently $F$, must be a biholomorphism. \qed
\medskip

\section{The proof of Theorem~\ref{T:short}}\label{S:short}

As $p$ is an orbit accumulation point, there is a point $a_0\in D_1$ and a sequence
$\{\phi_k\}\subset {\rm Aut}(D_1)$ such that $\phi_k(a_0)\to p$. Regardless of whether
$p$ is a peak point or $F(p)$ is a peak point, let us denote the relevant peak function as
$H$. Let $B$ denote a small ball centered at $p$, with $B\Subset U$, if $p$ is a peak point,
and centered at $F(p)$, with $B\Subset F(U)$, if $F(p)$ is a peak point. Depending on
whether $p$ or $F(p)$ is a peak point, set $G:=F^{-1}$ or $G:=F$, respectively. Finally,
set
\[
 h\,:=\,\begin{cases}
 		\left.H\circ G\right|_{\overline{B\cap D_2}}, &\text{if $p$ is a peak point,} \\
 		\left.H\circ G\right|_{\overline{B\cap D_1}}, &\text{if $F(p)$ is a peak point.}
 		\end{cases}
\]
By our hypothesis on $F$, it follows that $h$ satisfies all the conditions required 
of the function $h$ in Lemma~\ref{L:conv-const} for the appropriate choice of
$(D,p)$ depending on whether $p$ or $F(p)$ is a peak point.
\smallskip

Let us now denote the automorphisms discussed above as $\phi^1_k, \ k = 1,2,3,\dots$
Then, using $H$ or the function $h$ constructed above, depending on whether $p$
or $F(p)$ is a peak point, we deduce by
Lemma~\ref{L:conv-const} that $\phi_k^1 \lrarw \con{p}$ uniformly on compact
subsets of $D_1$.  Set $a_k := \phi^1_k(0)$. As $a_k\to p$, we may assume without
loss of generality that $a_k\in U$. Let $b_k := F(a_k)$, and let
$\phi^2_k\in {\rm Aut}(D_2)$ be an automorphism that maps $0$ to $b_k$ (which is
possible as ${\rm Aut}(D_2)$ acts transitively on $D_2$). Repeating the above argument,
$\phi_k^2 \lrarw \con{F(p)}$ uniformly on compact subsets of $D_2$. We have arrived
at the same result as in \eqref{E:CvgCon}. Thereafter, if we define
\[
  \OM_j^k := (\phi_k^j)^{-1}(\OM_j), \; \; j =1,2, \ k \in \Z_+,
\]
where $\OM_1 := U$ and $\OM_2 := F(U)$, then, reasoning exactly as in the
passage following \eqref{E:CvgCon}, we deduce that $(*)$ from Section~\ref{S:mThm}
holds true for our present set-up.
\smallskip

With $\{k_m\}_{m\in \Z_+}$ as given by $(*)$, let us define the maps $G_m : \OM^{k_m}_1\to
\OM^{k_m}_2$ by
\[
 G_m := (\phi_{k_m}^2)^{-1} \circ F \circ \phi_{k_m}^1. 
\]
By construction, each $G_m$ is a biholomorphic map. In particular,
\begin{equation}\label{E:atZero}
 G_m(0) = 0, \; \; \text{and} \; \; G'_m(0)\in GL(n,\C).
\end{equation}
We may assume, owing to \eqref{E:CvgCon}, that {\em the sequences $\{\OM^{k_m}_j\}_{m\in \Z_+}$
are increasing sequences}. By Montel's theorem, and arguing by induction, we can find sequences
$\{G_{l,m}\}$ and holomorphic maps $\Gamma_l : \OM^{k_l}_1\to \overline{D}_2$ such that:
\begin{itemize}
 \item $\{G_{1,m}\}_{m\in \Z_+}$ is a subsequence of $\{G_{\nu}\}_{\nu\in \Z_+}$ and
 $\{G_{l+1,m}\}_{m\in \Z_+}$ is a subsequence of $\{G_{l,\nu}\}_{\nu\in \Z_+}$;
 \item $\left.G_{l,m}\right|_{\OM^{k_l}_1}\lrarw \Gamma_l$, as $m\to \infty$, uniformly
 on compact subsets of $\OM^{k_l}_1$;
\end{itemize}
for $l = 1,2,3,\dots$ Owing to this construction, the rule
\[
 \Gamma(z)\,:=\,\Gamma_l(z) \; \; \text{if $z\in \OM^{k_l}_1$},
\]
gives a well-defined holomorphic map $\Gamma: D_1\to \overline{D}_2$.
\smallskip

Let us define $H_l := G_{l,l}$. Now suppose $\Gamma(D_1)\cap \bdy D_2\neq \varnothing$.
Then, $\exists\xi\in D_1$ such that $\Gamma(\xi)\in \bdy D_2$. Let $M\in \Z_+$ be so large
that $\OM^{k_M}_1 \ni \xi$. As $D_2$ is a bounded symmetric domain, it is taut. Thus, by
focusing attention on the sequence
\[
 \{\left.H_l\right|_{\OM^{k_M}_1} : l = M, M+1, M+2,\dots\} \subset \hol(\OM^{k_M}_1;D_2),
\]
we must conclude, by assumption, that $\Gamma(\OM^{k_M}_1)\subset \bdy D_2$. But, by
\eqref{E:atZero}, $\Gamma(0) = 0\notin \bdy D_2$. This is a contradiction, from which we infer:
\begin{itemize}
 \item[(a)] The range of $\Gamma$ is a subset of $D_2$.
\end{itemize}
Now observe that, by $(*)$, we have:
\begin{itemize}
 \item[(b)] The sequence $\{H_l: l = s, s+1, s+2,\dots\}$ converges uniformly to
 $\Gamma$ on $(1-1/s)\overline{D}_1$, $s\in \Z_+$.
 \item[(c)] $H_l^{-1}$ maps $0$ to $0$ and $(1-1/\nu^*(l))D_2$ into $D_1$ (since
 ${\sf dom}(H_l^{-1}) = {\sf range}(H_l) \supseteq \OM^{k_l}_2$).
\end{itemize}
In view of \eqref{E:atZero} and the fact that $D_1$ and $D_2$ are balanced, 
(a)-(c) are precisely the ingredients required to to repeat the argument in Step 1 of
the proof of Theorem~\ref{T:mThm} to infer that $\Gamma'(0)$ is invertible,
$\Gamma'(0)^{-1}: D_2\to D_1$ and
\[
 \Gamma'(0)^{-1}\circ \Gamma = {\sf id}_{D_1}.
\]

Thus, by (a), $\Gamma'(0)(D_1)\subset D_2$. One of the consequences of repeating the
argument contained in Step 1 in Section~\ref{T:mThm} is, in view of (c), that
$\Gamma'(0)^{-1}$ maps $(1-1/\nu^*(l))D_2$ into $D_2$ for every $l\in \Z_+$. As
$\nu^*$ is strictly increasing and $\Z_+$-valued, and as $\Gamma'(0)$ is injective,
this means that $\Gamma'(0)(D_1)$ contains $(1-1/s)D_2$ for {\em arbitrarily large}
$s\in \Z_+$, whence $\Gamma'(0)$ maps $D_1$ onto $D_2$. Hence,
$\Gamma'(0)|_{D_1}$ is a biholomorphism of $D_1$ onto $D_2$. \qed 
\medskip

\noindent {\bf Acknowledgements.} Jaikrishnan Janardhanan would like to thank his colleagues
and friends G.P. Balakumar, Dheeraj Kulkarni, Divakaran Divakaran and Pranav Haridas for many 
interesting discussions. Gautam Bharali contributed to this work while on sabbatical at the
Norwegian University of Science and Technology (NTNU) in Trondheim. He would like to acknowledge
the support and the hospitality of the Department of Mathematical Sciences at NTNU. He also thanks
John Erik Forn{\ae}ss for his helpful comments on this work.

\end{document}